\documentclass[12pt, reqno]{amsart}%
\usepackage[T1]{fontenc}
\usepackage{graphicx,subfigure}
\usepackage{epsfig}
\usepackage{amsmath}
\usepackage{amsfonts}
\usepackage{amssymb}
\usepackage{color}
\usepackage{enumerate}
\usepackage{rotating}
\usepackage{hyperref}
\usepackage{multirow}
\usepackage[overload]{empheq}

\providecommand{\U}[1]{\protect\rule{.1in}{.1in}}

\providecommand{\U}[1]{\protect\rule{.1in}{.1in}}
\numberwithin{figure}{section}
\numberwithin{table}{section} 
{\theoremstyle{plain}

} 
{\theoremstyle{definition}

}

\newcommand{\RR}{\mathbb{R}}
\newcommand{\RP}{\text{RP}}
\newcommand{\laL}{\lambda_{L}}
\newcommand{\laR}{\lambda_{R}}

\begin{document}
\title{Embedded delta shocks}

\address{\vspace{-18pt}
\newline Pablo Casta\~neda
\newline Department of Mathematics, ITAM
\newline R\a'io Hondo 1, Ciudad de M\'exico 01080, Mexico}
\email{pablo.castaneda@itam.mx}

\begin{abstract}
In 1977 Korchinski presented a new type of shock discontinuity in conservation laws. These singular solutions were coined $\delta$-shocks since there is a time dependent Dirac delta involved. A naive description is that such $\delta$-shock is of the overcompressive type: a two-family shock wave the four characteristic lines of which impinge into the shock itself. In this work, we open the fan of solutions by studying two-family waves without intermediate constant states but, possessing central rarefactions and also comprising $\delta$-shocks.
\end{abstract}

\date{Paper draft: October 2, 2019}

\subjclass[2000]{Primary: 35L65, 35L67; Secondary: 35M10, 58J45, 76T99.}
\keywords{Conservation laws, Riemann problem, $\delta$-shock, singular shocks, Rankine-Hugoniot, two-phase flow, chromatography.}

\maketitle

\section{Introduction}

\noindent
The introduction of $\delta$-shocks occurred forty years ago with the unpublished thesis \cite{Korchinski}, where such discontinuities appear in a theoretical context. Around that time, there was a simplified model for multiphase flow in porous media due to D.W. Peacement that also presented such a mass accumulation within one of these singularities, \cite{Dan}. Along these four decades, the applicability of $\delta$-shocks have emerged in many areas such as chromatography \cite{Key11,YaZh12}, magnetohydrodynamics \cite{NeOb08,TZZ94,WaYa18}, traffic flow \cite{LiXi18}, fluid dynamics \cite{KSZ04}, and perhaps also in flow in porous media \cite{BrMa07}, among other areas.

It is natural to consider a $\delta$-shock with speed $\sigma$ as an \textit{overcompressive shock wave}, which means a discontinuity satisfying that left and right characteristic lines impinge into the shock itself, \textit{i.e.},
\begin{equation}
 \lambda_{1,2}(U_L) \;>\; \sigma \;>\; \lambda_{1,2}(U_R),
 \label{eq:OC}
\end{equation}
for $U_L$ and $U_R$ the left and right Riemann data and $\lambda_{1,2}(U)$ the characteristic speeds for a point $U = (u,\,v)^T$ in state space; \textit{cf.} \cite{Key99,Korchinski,LiXi18,NeOb08,Sev07,TZZ94,ZhZh19}. Overcompressibility in Eq.~\eqref{eq:OC} is a natural extension of Lax classification, \cite{Lax57}, which considers also the following comparisons of speeds
\begin{alignat}{3}
 \lambda_{1,2}(U_L) \;>\; &\,\sigma \;>\; \lambda_{1}(U_R), & \qquad 
 \lambda_2(U_R) \;>\; &\,\sigma, 
 \label{eq:1Lax} \\
 \lambda_{2}(U_L) \;>\; &\,\sigma \;>\; \lambda_{1,2}(U_R), & \qquad &\,\sigma \;>\; \lambda_1(U_L) ,
 \label{eq:2Lax} \\
 \lambda_{2}(U_L) \;>\; &\,\sigma \;>\; \lambda_{1}(U_R), & \qquad
 \lambda_2(U_R) \;>\; &\,\sigma \;>\; \lambda_1(U_L).
 \label{eq:XS}
\end{alignat}
giving rise to \textit{1-Lax shock waves} in Eq.~\eqref{eq:1Lax}, \textit{2-Lax shock waves} in Eq.~\eqref{eq:2Lax} and, \textit{undercompressive} or \textit{transitional shock waves} in Eq.~\eqref{eq:XS}; for further details see \cite{Cas18,CAFM16} and references therein.

The types of shocks given by \eqref{eq:1Lax}-\eqref{eq:XS} are not found explicitly in the literature in conjunction to $\delta$-shocks.
From extensive large bibliographic review in \cite{YaZh12} for models with $\delta$-shocks, we notice that the analyzed and identified conservation laws models are weakly coupled and of the form
\begin{equation}
 \begin{array}{rclrclr}
  u_t + \big(F(u,\,v)\big)_x &=& 0, &
  (u^\alpha v)_t + \big(G(u,\,v)\big)_x &=& 0, & (x,\,t) \in \RR\times\RR^+,
 \end{array}
 \label{sys:linear} 
\end{equation}
where $\alpha \in \{0,\,1\}$ and, $F$ and $G$ are linear in $v$, see also \cite{Key11}.

Consider the case $\alpha = 0$ and notice that for a Riemann problem including a $\delta$-shock, the shock speed is extracted from (\ref{sys:linear}.a),  which determines left and right transport speeds $c_L = G(u_L,\,v)/v$ and $c_R = G(u_R,\,v)/v$ for Eq.~(\ref{sys:linear}.b). 
Now, an equation of the transport type $v_t + cv_x = 0$ should be solved at left and right of $x = \sigma t$, with $c = c_L$ and $c_R$, respectively. The characteristic lines from (\ref{sys:linear}.a) impinge into the shock wave, however, comparisons of $\sigma$ against $c_L$ and $c_R$ are free, and then inequalities \eqref{eq:OC}-\eqref{eq:XS} may hold; necessarily the compressibility is preserved. For $v$ we have two transport equations, which can only carry information from the Riemann data; the $\delta$-shock is consequence solely of the imbalance of mass at $x = \sigma t$. 
Still, this $\delta$-shock is surrounded by constant states rather than rarefaction waves.

An overcompressive shock is a restrictive wave in the sense that it is an isolated discontinuity for a Riemann problem connecting left and right states $U_L$, $U_R$ via this shock; Eq.~\eqref{eq:OC} holds, and there can be neither preceding nor succeeding waves, only constant states on both sides of the discontinuity.
Our main result is the construction of the other types of shock waves related to \eqref{eq:1Lax}-\eqref{eq:XS} with a $\delta$-shock involved.
The new $\delta$-shocks may precede or succeed rarefaction waves. Hence, classical Riemann solutions with two wave groups. Typically, there exists an intermediate constant state separating wave groups. The authors in \cite{SiMa14} took the endeavor to produce a set of conservation law models possessing Riemann solutions without such intermediate constant states. Remarkably, the solutions we present here possess a $\delta$-shock rather than these intermediate constant states. Another directions are given in \cite{Che19}, where Riemann solutions are reported that possess no intermediate constant states but $\delta$-contact discontinuities and, in \cite{NeOb08}, where  interaction of classical waves and $\delta$-shocks is given in a positive time.

The rest of this work is organized as follows.
In Sec.~\ref{sec:Kor}, we reconstruct the overcompressive shock wave found by Korchinski. In Sec.~\ref{sec:Rars}, we present the new $\delta$-shocks of type \eqref{eq:1Lax}-\eqref{eq:XS} with preceding or succeeding central rarefaction fans. Finally, in Sec.~\ref{sec:example}, we present a Riemann solution possessing two $\delta$-shocks. 
Some concluding remarks are presented in Sec.~\ref{sec:conclusion}.

\subsection{The first analysis, back to 1977}
\label{sec:Kor}

\noindent
Take Korchinski system \cite{Korchinski}, and rescale it as in \cite{Key99}:
\begin{equation}
 \begin{array}{rclrclr}
  u_t + (u^2)_x &=& 0, & 
  v_t + (uv)_x  &=& 0, & (x,\,t) \in \RR\times\RR^+.
 \end{array}
 \label{sys:Kor}
\end{equation}
We denote Riemann problems as $\RP(U_L,\,U_R)$, comprising a system of conservation laws (as \eqref{sys:Kor}, \eqref{sys:mod} or \eqref{sys:2deltas}), and a discontinuous initial condition
\begin{equation}
 U(x,\,0) \;=\;
 \begin{cases}
  U_L, & x < 0, \\
  U_R, & x > 0.
 \end{cases}
 \label{eq:IC}
\end{equation}
From the well-known Rankine-Hugoniot condition, a shock front for $u$ with propagation speed $\sigma = u_L + u_R$ exists when $u_R < u_L$ holds. In the presence of this shock wave, $v$ changes across the front line. The solution profile can be written as
\begin{equation}
  v(x,\,t) = v_L + (v_R - v_L)\mathcal{H}(x - \sigma t) + k(t)\delta(x - \sigma t),
  \label{sol:Kor}
\end{equation}
where $\mathcal{H}$ is the Heaviside step function and $\delta$ is the Dirac delta, see \cite{Key99,Korchinski}.

In a conservation law, the change of mass in an interval is equal to the net flow of mass at the boundary. For an interval $x \in [a,\,b]$ with $a \ll 0 \ll b$, the mass balance of $v(x,\,t)$ in \eqref{sol:Kor} is given by
\begin{eqnarray}
u_L v_L - u_R v_R
  &=& \frac{d}{dt}\int_{a}^{b} v(x,\,t)\,dx
  \;=\; \frac{d}{dt} \left[ \int_{a}^{\sigma t} v_L\,dx
     + \int_{\sigma t}^b v_R\,dx
     + \int_a^b k(t)\delta(x - \sigma t) \,dx \right]
     \nonumber \\
  &=& \sigma(v_L - v_R) + k'(t).
  \label{eq:balanceKor}
\end{eqnarray}
Equating these equalities and integrating over $t$ leads to
$k(t) = (u_Rv_L - u_Lv_R)t$, since the initial condition \eqref{eq:IC} implies $k(0) = 0$.
Thus, this Riemann problem has solution
$$U(x,\,t) \;=\; \left(
\begin{array}{c}
 u_L + (u_R - u_L)\mathcal{H}(x - \sigma t) \\
 v_L + (v_R - v_L)\mathcal{H}(x - \sigma t) + (u_Rv_L - u_Lv_R)t\,\delta(x - \sigma t)
\end{array}
\right),$$
which is plotted in Fig.~\ref{fig:profiles}. The second coordinate state possesses a $\delta$-shock with growing amplitude $k(t)$. 

Of course, these computations hold in the sense of distributions, see \cite{Che19,SaSe02,She18}. However, the Riemann solutions in the following sections comprise rarefactions that are difficult to handle in these distributions. Even if it is possible to compute the generalized Rankine-Hugoniot conditions given in \cite{TZZ94}, see also \cite{Che19}, we prefer for simplicity direct computations as in \eqref{eq:balanceKor}.

\section{A $\delta$-shock near a rarefaction wave}
\label{sec:Rars}

\noindent
In this section we modify system \eqref{sys:Kor} in order to produce a richer set of discontinuities around a $\delta$-shock. We consider
\begin{equation}
 \begin{array}{rclrclr}
  u_t + (u^2)_x  &=& 0, & 
  v_t + (uv^2)_x &=& 0, & (x,\,t) \in \RR\times\RR^+.
 \end{array}
 \label{sys:mod}
\end{equation}
As before, from (\ref{sys:mod}.a), a solution for the $\RP(U_L,\,U_R)$ has a shock wave with speed $\sigma = u_L + u_R$ when $u_R < u_L$; this fact will be assumed from now on.

Now, the nonlinear flux for $v$ is $uv^2$, so at constant $U_{L,R}$ we have characteristic speeds, $\laL = 2u_Lv_L$ at the left of the shock front and $\laR = 2u_Rv_R$ at the right. (The other two characteristic speeds satisfy $\tilde{\lambda}_L := 2u_L > \sigma > 2u_R =: \tilde{\lambda}_R$.) 
In the original model, the flux for $v$ is linear around the shock and the $\delta$-shock is a consequence of this imposed transport.

\begin{figure}[htb]
  \begin{picture}(480,200)
    \put(15,0){\resizebox{0.94\textwidth}{!}{\includegraphics{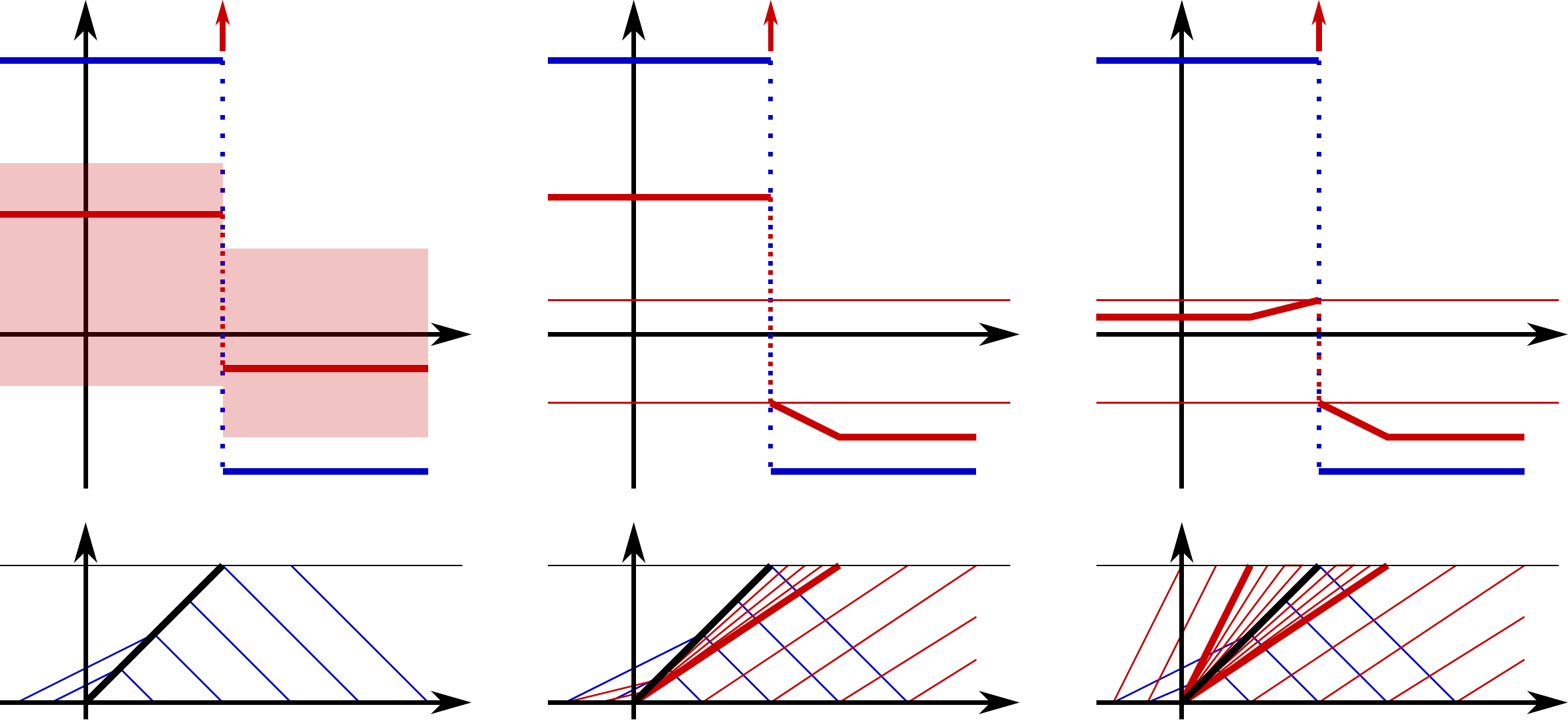}}}
    \put( 45,195){$u,\,v$}
    \put( 82,191){$\delta$}
    \put(  0,184){$u_L$}
    \put(  0,139){$v_L$}
    \put(138, 67){$u_R$}
    \put(138, 96){$v_R$}
    \put(150,105){$x$}
    \put( 80, 47){$\sigma$}
    \put( 45, 47){$t$}
    \put(150,  2){$x$}
    \put(200,195){$u,\,v$}
    \put(237,191){$\delta$}
    \put(155,184){$u_L$}
    \put(155,145){$v_L$}
    \put(293, 67){$u_R$}
    \put(293, 77){$v_R$}
    \put(304,105){$x$}
    \put(233, 47){$\sigma$}
    \put(248, 47){$\lambda_R$}
    \put(200, 47){$t$}
    \put(304,  2){$x$}
    \put(353,195){$u,\,v$}
    \put(391,191){$\delta$}
    \put(310,184){$u_L$}
    \put(310,112){$v_L$}
    \put(447, 67){$u_R$}
    \put(447, 77){$v_R$}
    \put(459,105){$x$}
    \put(364, 47){$\lambda_L$}
    \put(387, 47){$\sigma$}
    \put(402, 47){$\lambda_R$}
    \put(353, 47){$t$}
    \put(459,  2){$x$}
  \end{picture}
  \caption{\small Profiles with $\delta$-shocks. We use blue, red and black for curves related to $u$, $v$, and both $u$ and $v$. On top, solid lines represent constant states and rarefactions, dotted lines are shock waves at $x = \sigma t$ (arrows with $\delta$ are schematic directions of $\delta$-shocks); on bottom, we have characteristic speeds on $xt$ plane, the horizontal thin line is time $t = 1$ taken as reference for the advance of waves on top panel profiles; $x = \sigma t$ is in thick dark line, $x = \lambda_{\{L,\,R\}}t$ are in thick red lines.
  All Riemann problems have $u_L > u_R$. Left panels is RP for \eqref{sys:Kor}, shaded regions represent that this configuration exists for any choice of $v_L$, $v_R$. Central and right panels are RP for \eqref{sys:mod}, thin horizontal red lines represent the thresholds $\lambda_L = \sigma$ and $\lambda_R = \sigma$; $\lambda_L < \sigma$ implies a rarefaction before the $\delta$-shock as in right panels, similarly $\lambda_R > \sigma$ implies rarefaction after the $\delta$-shock as in central and right panels.}
  \label{fig:profiles}
\end{figure}
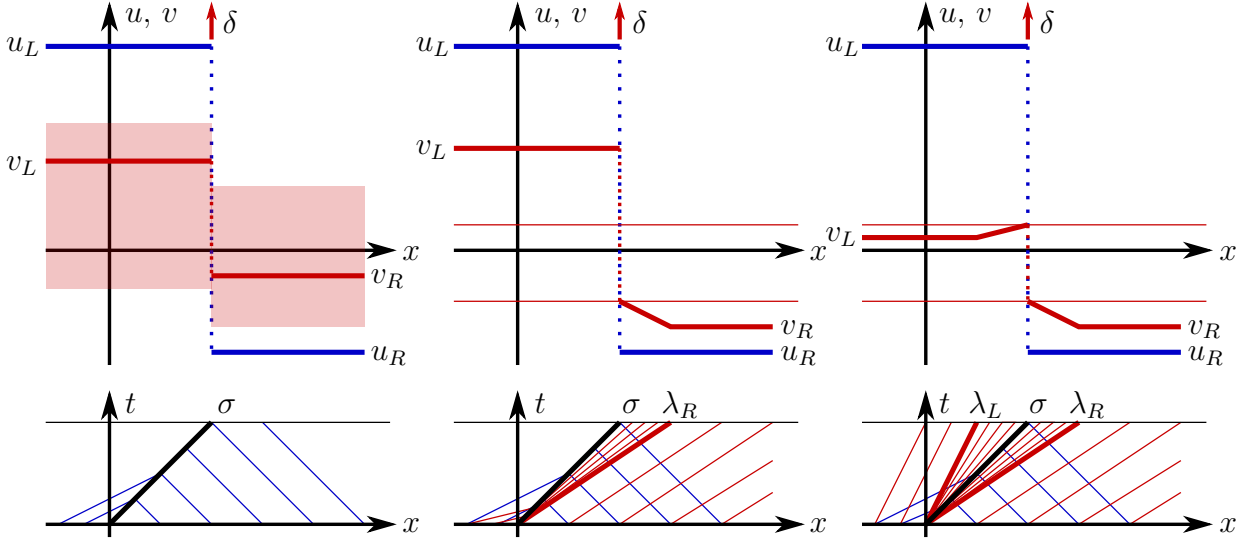

New scenarios arise when $\laL,\,\laR > \sigma$ as in Eq.~\eqref{eq:1Lax}, $\sigma > \laL,\,\laR$ as in Eq.~\eqref{eq:2Lax}, or $\laL < \sigma < \laR$ as in Eq.~\eqref{eq:XS}. We study the first and third cases; the second case is similar to the first one. Notice that in the first case, as $\sigma < \lambda_R$, the gap in characteristic lines in $xt$ plane can be filled with a centered rarefaction fan via the nonlinear flux in (\ref{sys:mod}.b). In the third case $\lambda_L < \sigma < \lambda_R$ hold, thus preceding and subsequent rarefactions appear around the $\delta$-shock, see bottom panels in Fig.~\ref{fig:profiles}.

\subsection{The case of ${\delta}$-shock\,--\,rarefaction}

\noindent
When the speed inequalities $\sigma < \laL,\,\laR$ hold, at the left of the shock discontinuity, the result must be as in the Korchinski case: $\laL,\,\tilde{\lambda}_L > \sigma$. However, at the right of this shock a rarefaction must appear to fill the gap between $\sigma t$ and $\laR t$ in $xt$ plane. For this reason, we take the solution \textit{ansatz}
\begin{equation}
  v(x,\,t) 
   \;=\; v_L 
       + \left(\frac{x/t}{2 u_R} - v_L\right)\mathcal{H}(x - \sigma t)
       + \left(v_R - \frac{x/t}{2 u_R}\right)\mathcal{H}(x - \lambda_R t)
       + k(t)\delta(x - \sigma t),
  \label{sol:mod1}
\end{equation}
comprising a ``fast'' rarefaction that also satisfies (\ref{sys:mod}.b).
As in \eqref{eq:balanceKor}, the mass balance is computed from \eqref{sol:mod1} as
\begin{eqnarray*}
u_L v_L^2 - u_R v_R^2 
   &=& \frac{d}{dt} \left[ \int_{a}^{\sigma t} v_L\,dx
       + \int_{\sigma t}^{\lambda_R t} \frac{x/t}{2u_R}\,dx
       + \int_{\lambda_R t}^b v_R\,dx
       + \int_a^b k(t)\delta(x - \sigma t) \,dx \right] \nonumber \\
  &=& \sigma v_L
      + \frac{\lambda_R^2 -  \sigma^2}{4u_R}
      - \lambda_Rv_R + k'(t),
\end{eqnarray*}
which leads to $k(t) = \left[u_Lv_L^2 - \sigma v_L + \sigma^2/(4u_R)\right]t.$
An example with $U_L = (2,\,1)^T$, $U_R = (-1,\,-3/4)^T$ is given in Fig.~\ref{fig:profiles}.

Borrowing terminology from Riemann problems for conservation laws (see \cite{Lax57,Liu74}), we say that this solution is given by a $\delta$-shock of type 1-Lax for the first wave group (\emph{i.e.}, characteristic speeds satisfy \eqref{eq:1Lax}), the second wave group is a second family (or fast) rarefaction. This 1-Lax $\delta$-shock possesses a linearly increasing Dirac delta, as the one in the Korchinski model, see Eq.~\eqref{eq:balanceKor}. Moreover, notice the lack of intermediate constant state between wave groups.

\subsection{The case of rarefaction\,--\,${\delta}$-shock\,--\,rarefaction}
\label{sec:deltaX}

\noindent
We consider now the case $\laL < \sigma < \laR$. The \textit{ansatz} satisfying (\ref{sys:mod}.b) is
\begin{eqnarray}
  v(x,\,t)
  &=& v_L 
       + \left(\frac{x/t}{2 u_L} - v_L\right)\mathcal{H}(x - \lambda_L t)
       + \left(\frac{x/t}{2 u_R} - \frac{x/t}{2 u_L}\right)\mathcal{H}(x - \sigma t) \nonumber \\
  & &  + \left(v_R - \frac{x/t}{2 u_R}\right)\mathcal{H}(x - \lambda_R t)
       + k(t)\delta(x - \sigma t),
  \label{sol:mod2}
\end{eqnarray}
which comprises ``slow'' and ``fast'' rarefactions.
The mass balance is computed from \eqref{sol:mod2} as
\begin{eqnarray*}
u_L v_L^2 - u_R v_R^2 
   &=& \frac{d}{dt} \left[ \int_{a}^{\lambda_L t} \!\!\!\!\! v_L\,dx
       + \! \int_{\lambda_L t}^{\sigma t} \! \frac{x/t}{2u_L}\,dx
       + \! \int_{\sigma t}^{\lambda_R t} \! \frac{x/t}{2u_R}\,dx
       + \! \int_{\lambda_R t}^b \!\!\!\! v_R\,dx
       + \! \int_a^b \!\!\! k(t)\delta(x - \sigma t) \,dx \right] \nonumber \\
   &=& \lambda_L v_L
      + \frac{\sigma^2 - \lambda_L^2}{4u_L}
      + \frac{\lambda_R^2 - \sigma^2}{4u_R}
      - \lambda_Rv_R + k'(t) 
      \nonumber \\
\end{eqnarray*}
which leads to $k(t) = \sigma^2(u_L - u_R)/(4 u_L u_R)t.$
Notice that stationary shocks, \textit{i.e.} with $\sigma = u_L + u_R = 0$, do not produce deltas.
An example with $\sigma = 1$: $U_L = (2,\,1/8)^T$, $U_R = (-1,\,-3/4)^T$ is given in Fig.~\ref{fig:profiles}.

This solution is given by a first family (or slow) rarefaction as first wave group, a $\delta$-shock of transitional type, see \eqref{eq:XS}, and a second family (or fast) rarefaction as the second wave group. Notice the linear behaviour of $k(t)$ and the lack of intermediate constant state between wave groups.

\section{Example of a wave with two $\delta$-shocks}
\label{sec:example}

\noindent
In previous sections we have studied wave groups possessing a single $\delta$-shock. Our aim now is to construct a new model supporting two of such singular discontinuities. This model possesses the features of models in \cite{YaZh12}.

\begin{figure}[htb]
\vspace{-0.3cm}
\resizebox{0.55\textwidth}{!}{
  \begin{picture}(350,80)
    \put(0,0){\resizebox{0.85\textwidth}{!}{\includegraphics{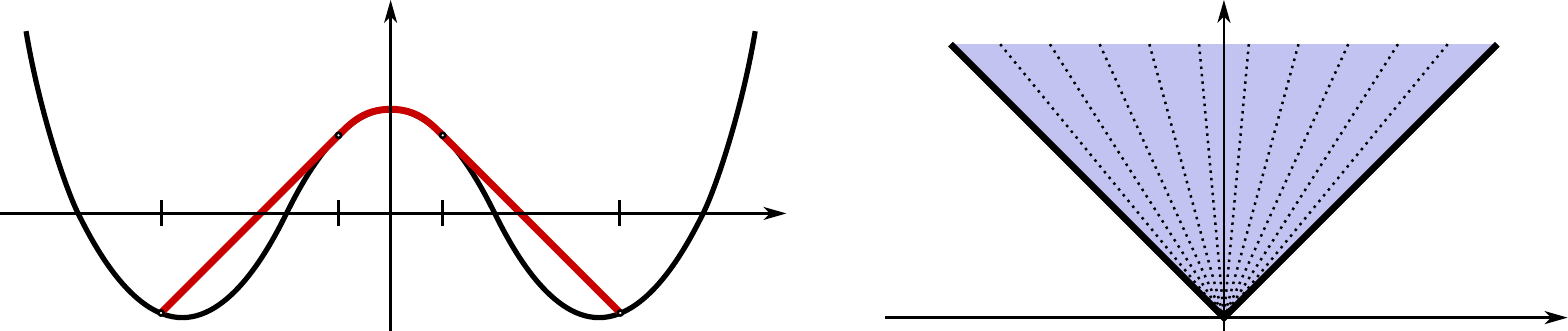}}}
    \put(193,17){\large $u$}
    \put( 72,75){\large $f(u)$}
    \put( 35,17){\large $u_R$}
    \put( 80,17){\large $a_R$}
    \put(108,17){\large $a_L$}
    \put(153,17){\large $u_L$}
    \put(400, 0){\large $x$}
    \put(300,75){\large $t$}
    \put(230,32){\large $x = \sigma^- t$}
    \put(350,32){\large $x = \sigma^+ t$}
  \end{picture}
}
\hfill
\resizebox{0.34\textwidth}{!}{
  \begin{picture}(400,160)
    \put(0,0){\resizebox{0.85\textwidth}{!}{\includegraphics{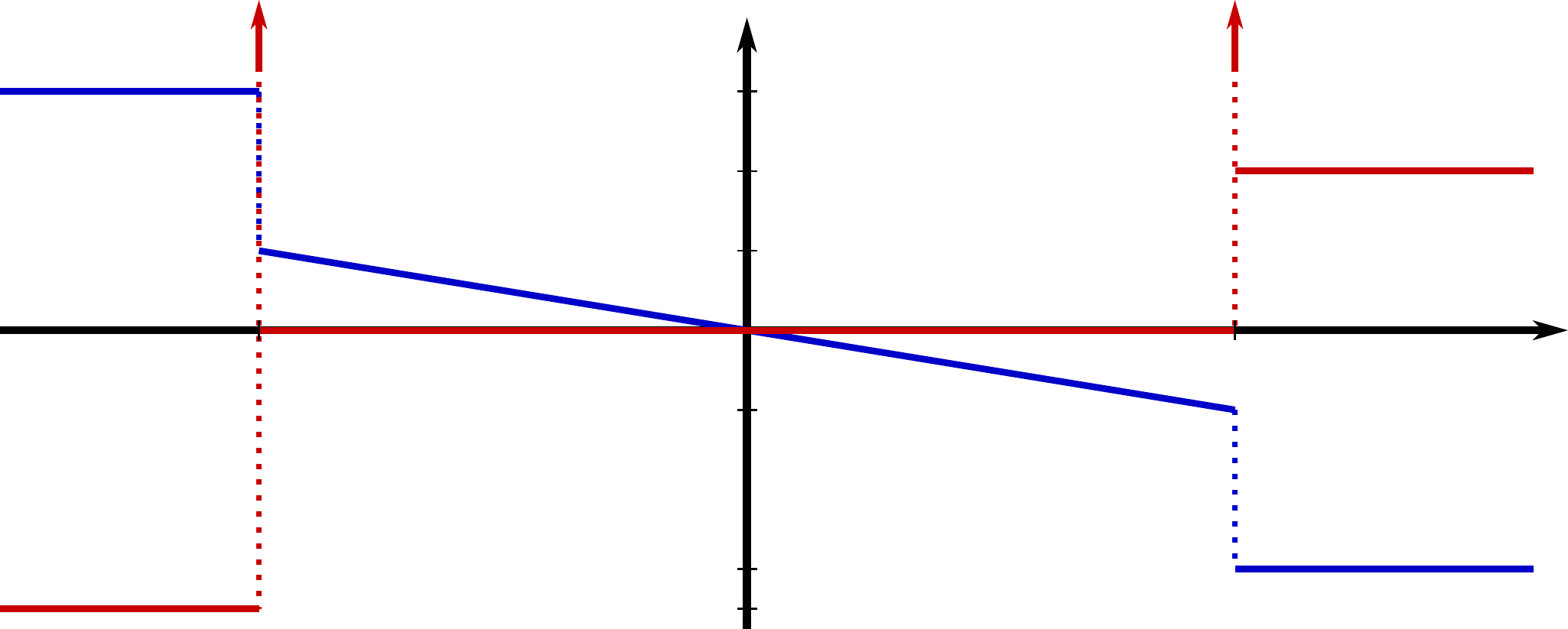}}}
    \put(385, 58){\Huge $x$}
    \put(195,140){\Huge $u,\,v$}
    \put(18,145){\Huge $u_L$}
    \put(18, 13){\Huge $v_L$}
    \put(340, 23){\Huge $u_R$}
    \put(340,125){\Huge $v_R$}
    \put( 72,142){\Huge $\delta^-$}
    \put(320,142){\Huge $\delta^+$}
  \end{picture}
}
  \caption{\small Left: Flux function \eqref{eq:flux} in black, Ole\u{\i}nik convex hull for $u_L = 1 = -u_R$ in red; the envelope is tangent at $a_L$ and $a_R$. 
  Center: Characteristic speeds for the associated RP, solid lines represent shock waves, dotted lines represent centered rarefaction fan.
  Right: Profile solution for system~\eqref{sys:2deltas}; blue is $u(x,\,t)$ profile, red is $v(x,\,t)$ profile. Two $\delta$-shocks at $\sigma^\pm t$, the ``amplitudes'' are specified; $\delta^\pm$ denotes the pulse $k_\pm(t)\delta(x - \sigma^\pm t)$.}
  \label{fig:flux}
\end{figure}
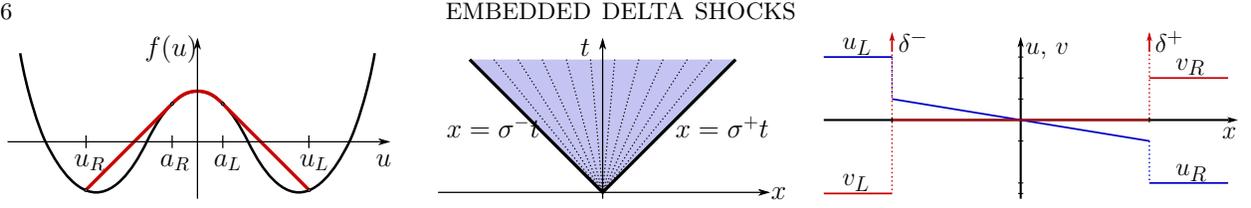

Let us take a modification of \eqref{sys:Kor} with a distinguished conservation for $u$ and repeat the conservation law for $v$, see (\ref{sys:Kor}.b). We write the system
\begin{equation}
 \begin{array}{rclrclr}
  u_t + f(u)_x &=& 0, & 
  v_t + (uv)_x  &=& 0, & (x,\,t) \in \RR\times\RR^+,
 \end{array}
 \label{sys:2deltas}
\end{equation}
where the flux $f(u)$ is a double-well function. For the sake of simplicity, from here and on, we consider
\begin{equation}
 f(u) \;=\;
  \begin{cases}
   (u + 2)^2 - 1, &\text{for }\; u < -1 \\
   u^2 + 1,       &\text{for }\; u \in [-1,\,1] \\
   (u - 2)^2 - 1, &\text{for }\; u > 1
  \end{cases},
 \label{eq:flux}
\end{equation}
and for the Riemann problem, we consider $u_L = - u_R = (3 + \sqrt{2})/2$. Then, the solution for $u$ is
\begin{equation}
 u(x,\,t) \;=\;
 u_L
 + \left(\frac{x/t}{a_L} - u_L\right)\mathcal{H}(x - \sigma^-t)
 + \left(u_R - \frac{x/t}{a_L}\right)\mathcal{H}(x - \sigma^+t)
 \label{eq:usol}
\end{equation}
where from Ole\u{\i}nik construction ({\it cf.} \cite{Cas16} and Fig.~\ref{fig:flux}), we have $\sigma^+ = -\sigma^- = 1$ and $a_R = -a_L = 1/2$.

The flux for $v$ is $u$, thus from \eqref{eq:usol} we notice that such flux is zero at $x = 0$.
The \textit{ansatz} for this system is
\begin{eqnarray}
  v(x,\,t) &=& 
  v_L + (0 - v_L)\mathcal{H}(x - \sigma^-t) + (v_R - 0)\mathcal{H}(x - \sigma^+t) \nonumber \\
  & & + k_-(t)\delta(x - \sigma^- t)
  + k_+(t)\delta(x - \sigma^+ t),
  \label{sol:mod3}
\end{eqnarray}
the solution of which fulfills (\ref{sys:2deltas}.b) and \eqref{eq:usol}. Indeed, the constant regions for $x \notin [\sigma^- t,\, \sigma^+ t]$ satisfy directly $v_t = (uv)_x = 0$. For $x \in (\sigma^- t,\, \sigma^+ t)$, we have from \eqref{eq:usol} and assuming it must be a rarefaction, that it has the form $v(x,\,t) = mx/t$ for a slope $m$ to be specified. Then, by substituting this form into (\ref{sys:2deltas}.b) we obtain
$$v_t + (u(x,\,t)\,v)_x \;=\; -\frac{mx}{t^2} + \frac{1}{a_Lt} \frac{mx}{t} + \frac{x}{a_Lt} \frac{m}{t} \;=\; \frac{mx}{t^2}\left(-1 + \frac{2}{a_L}\right) \;=\; 0,$$
which holds only for $m = 0$.

Considering the positive axis, the change of mass of $v$ for $x \geq 0$ is given from \eqref{sol:mod3} as
\begin{equation*}
0 - u_L v_L 
  \;=\; \frac{d}{dt} \left[ \int_0^{\sigma^+ t} 0\,dx
        + \int_{\sigma^+ t}^b v_R\,dx 
        + \int_0^b k_+(t)\delta(x - \sigma^+ t) \,dx \right]
  \;=\; -\sigma^+ v_R + k'_+(t).
\end{equation*}
Thus, $k_+(t) = (\sigma^+ v_R - u_R)t$, and $k_-(t) \;=\; -(\sigma^- v_L - u_L)t,$ from an analogous treatment for the change of mass of $v$ for $x \leq 0$.

In Fig.~\ref{fig:flux} we plot the solution profile for $\RP(U_L,\, U_R)$, where $U_L = (-(3+\sqrt{2})/2,\, v_L)^T$ and $U_R = ((3+\sqrt{2})/2,\, v_R)^T$, for $v_L < u_R$ and $v_R < u_L$; for these settings $k_+(t), k_-(t) > 0$ for all times, the amplitude of both $\delta$-shocks is positive.

\section{Concluding remarks}
\label{sec:conclusion}

\noindent
A crucial feature in constructing the solutions in Sec.~\ref{sec:Rars} is the nonlinear behaviour of $G(u,\,v)$ in $v$, see \eqref{sys:linear}. From Eq.~(\ref{sys:mod}.a), or similar, we can extract the speed $\sigma$, which determines the existence and localization of $\delta$-shocks. The second flux, \emph{i.e.} $G(u,\,v)$, establishes thresholds by comparing $\lambda_L = G_v(u_L,\,v_L)$ and $\lambda_R = G_v(u_R,\,v_R)$ to $\sigma$. Notice that $v(x,\,t) \to \lambda_L$ ($\lambda_R$, respectively) as $x \to \sigma t-$ ($\sigma t+$, resp.), so a ``transitional $\delta$-shock'' has zero amplitude when $\lambda_L = \lambda_R = \sigma$ hold, but there is a bump at $x = 0$ (typically $v(0,\,t) = 0$ is larger than $v_L,\,v_R$). 
In such a situation a $\delta$-shock is masked within a bump; small perturbations of the Riemann data will reproduce the linear growing of the delta.
In other words, $\delta$-shocks can be masked with specific mathematical settings, which stands in contradistinction to their nature from the physical point of view, this reinforces the idea of possible $\delta$-shocks not reported in the literature.

In \cite{LeF90}, LeFloch established the existence of solutions for Cauchy problems in a model similar to \eqref{sys:2deltas} for convex flux $f(u)$. For such fluxes, the Riemann problem may possess a single $\delta$-shock. Here we constructed an elegant solution comprising two $\delta$-shocks. In \cite{Che19}, a solution with three $\delta$-shocks appears for a $3 \times 3$ system of conservation laws. Actually, following the ideas in Sec.~\ref{sec:example}, we can give a flux $f(u)$ that allows the generation of any number of $\delta$-shocks; each contact discontinuity from the Ole\u{\i}nik convex hull construction may become a $\delta$-shock.

On the other hand, solutions comprising rarefactions and $\delta$-shocks were presented in Sec.~\ref{sec:Rars}, and we noticed the absence of intermediate constant states in all of them. In the Riemann solutions foreseen in classical theory  by Lax and Liu (\emph{cf.} \cite{Lax57,Liu74}), the existence of intermediate constant states is necessary for the structural stability. The lack of these states is rather covered by $\delta$-shocks which suggest to be more general than the alternative of the ``organizing center'' given in \cite{SiMa14}. 
In summary, we can construct $2 \times 2$ Riemann solutions with any number of $\delta$-shocks, we can further compose them with different waves, as the transitional shocks used in \cite{Lozano}. Therefore, we have an eye-catching phenomenon that emerges with potential giving use to new solutions. These solutions arise in stark contrast of what is know for strictly hyperbolic systems of conservation laws.



\end{document}